\documentclass[12pt]{article}

\usepackage{amsmath}
\usepackage{amsfonts}
\usepackage{amssymb}
\usepackage{amscd}
\usepackage{latexsym}
\usepackage{mathrsfs}
\usepackage{amsthm}

\newtheorem{Corollary}{Corollary}
\newtheorem{Theorem}{Theorem}
\newtheorem{Proposition}{Proposition}
\newtheorem{Definition}{Definition}

\def\BB{\mathscr{B}}

\def\de{\mathrm d}
\def\D{\mathcal{D}}
\def\e{\varepsilon}
\def\E{\mathcal{E}}
\def\F{\mathcal{F}}
\def\G{\mathcal{G}}
\def\H{\mathcal{H}}
\def\M{\mathcal{M}}
\def\N{\mathcal{N}}
\def\R{\mathbb{R}}
\def\S{\mathcal{S}}
\def\T{\mathcal{T}}
\def\im{\mathop{\mathsf{Im}}\nolimits}
\def\re{\mathop{\mathsf{Re}}\nolimits}

\begin{document}

\title{Some improvements in the method \\ of the weakly conjugate operator}

\author{Serge Richard}
  \date{\small
    \begin{quote}
      \emph{
    \begin{itemize}
    \item[]
			Institut Camille Jordan,
			Universit\'e Claude Bernard Lyon 1, \\
			43 avenue du 11 novembre 1918,
			69622 Villeurbanne cedex, France
    \item[]
      \emph{E-mails\:\!:}
      srichard@math.univ-lyon1.fr
    \end{itemize}
      }
    \end{quote}
    October 2005
  }

\maketitle

\begin{abstract}
We present some improvements in the method of the weakly conjugate operator,
one variant of the Mourre theory.
When applied to certain two-body Schr\"odinger operators, this leads
to a limiting absorption principle that is uniform on the positive real axis.
\end{abstract}

\noindent
{\bf Mathematics Subject Classification (2000):} 35P05, 46L60, 47A10, 81Q10.

\noindent
{\bf Key words:} conjugate operator methods, limiting absorption principle,
zero energy, two-body Schr\"odinger operator.

\section{Introduction}
Recently there has been an increasing interest in the study of two-body Schr\"odinger 
operators near the threshold at energy zero (see for example \cite{FS} or \cite{JN}).
Since a positive commutator in the sense of Mourre does not exist
at this energy, the usual method of the conjugate operator can not be used 
in that particular situation.
On the other hand, the method of the weakly conjugate operator gives the 
existence of the boundary values of resolvents also at thresholds but 
applies only to situations where the operators have no bound states at all.
However the authors of \cite{FS} derive a limiting absorption principle at zero energy
for a special class of two-body Schr\"odinger operators which have bound states
below zero.
In this context an improvement of the method of the weakly conjugate operator
that will cover the behaviour at thresholds of operators with bound states would
be of interest. The purpose of this letter is to describe such an extension
and to give an application to two-body Schr\"odinger operators.

Let us recall the main idea of methods based on a conjugate operator. 
One way to obtain strong results for the spectral analysis of a self-adjoint
operator $H$ is to find an auxiliary self-adjoint operator $A$
such that the commutator $[iH,A]$ is positive in a suitable sense.
In the method of the conjugate operator one looks for intervals $J$ of $\R$ 
such that 
\begin{equation}\label{Mou}
E(J)[iH,A]E(J)\ \geq \ a\;\!E(J)
\end{equation} 
for some strictly positive constant $a$ that depends on $J$, where $E(J)$ 
denotes the spectral projection of $H$ on the interval $J$.
For the method of the weakly conjugate operator one assumes that $[iH,A]>0$,
{\it i.e.}~the commutator is positive and injective.
This requirement is closer to the initial Kato-Putnam theory, on which 
it improves.

The first approach has reached a very high degree of precision 
and abstraction in \cite{ABG}. 
There also exists a huge number of applications based on an inequality of the 
form \eqref{Mou}.
The second approach was initiated in \cite{BKM} and fully developed 
in \cite{BM}. Only a few papers contain applications, 
see for example \cite{IM}, \cite{MP} or \cite{MR}.
We also mention \cite{FS} and \cite{H} that contain arguments that are very 
close to this method.
Its main disadvantage is that if the method can be applied
to $H$, then the spectrum of $H$ is purely absolutely continuous, 
which limits drastically its range of applications.
On the other hand, it leads to a limiting absorption principle that
is uniform on $\R$ and to global $H$-smooth operators, that are of special 
interest.
We refer to \cite{S} and references therein for more information on that subject.

Motivated by some calculations borrowed from \cite{FS} we shall
prove in this letter that the fundamental assumption $[iH,A]>0$ of the method of
the weakly conjugate operator can be weakened. 
The main idea is that $H$ itself can add some positivity.
Surprisingly, the new requirement is that there exists a constant $c\geq 0$
such that 
\begin{equation*}
-c\;\!H + [iH,A] \ > \ 0\ .
\end{equation*}
This inequality together with some technical assumptions lead
to a limiting absorption principle that is either uniform on $\R$ if
$c=0$ or uniform on $[0,\infty)$ if $c>0$.
The absolute continuity of the spectrum and $H$-smooth operators
are then standard byproducts of that estimate.

In the next section we introduce the framework and state the abstract 
result. Its proof is postponed until Section \ref{secproof}.
In-between, we give an application to two-body Schr\"odinger operators.
We prove that under suitable conditions such operators admit a limiting
absorption principle uniform on $[0,\infty)$. 
Since our approach applies to operators that may have discrete spectrum 
below zero, our abstract result is really an improvement 
of the method developed in \cite{BKM} and \cite{BM}.  

We close this introduction with some comments on generality. 
As in the early papers on the method of the conjugate
operator, our condition on the second commutator $\big[i[iH,A],A\big]$
can certainly be weakened. Also an approach divided into two stages
(first by dealing with bounded operators and then by applying the result
to the resolvent $(H-\lambda_0)^{-1}$ for a real $\lambda_0$ outside 
the spectrum of $H$) would certainly lead to some improvements.
However, since such modifications would also lengthen and complicate our arguments, 
we decided not to take them into account in this letter.

\section{The abstract construction}\label{secintro}

Let $\H$ be a Hilbert space with scalar product $\langle \cdot, 
\cdot \rangle$ and norm  $\| \cdot \|$. 
We consider a self-adjoint operator $H$ in $\H$ with its domain denoted by
$\G^2$ and its form domain denoted by $\G^1$.
Endowed with the corresponding graph norms, $\G^2$ and $\G^1$ 
are also Hilbert spaces. 
Their adjoint spaces (topological anti-duals) are denoted by $\G^{-2}$ and 
$\G^{-1}$, and by identifying $\H$ with its adjoint through the Riesz isomorphism 
one has the continuous dense embeddings~:
\begin{equation*}
\G^2 \hookrightarrow \G^1 \hookrightarrow \H \hookrightarrow \G^{-1} 
\hookrightarrow \G^{-2}.
\end{equation*}

Let $\{W_t\}_{t \in \R}$ be a strongly continuous unitary group in $\H$ with 
its self-adjoint generator denoted by $A$.
We assume that for each $t \in \R$, $W_t$ leaves $\G^2$ invariant.
It is then a standard fact that $\{W_t\}_{t \in \R}$ induces a $C_0$-group in 
each space $\G^s$ introduced above \cite[Sec.~6.3]{ABG}.
We keep the same notation for these $C_0$-groups.

Now, let us consider an operator $S \in \BB(\G^1, \G^{-1})$ that satisfies
$S>0$, {\it i.e.}~$\langle f,Sf \rangle>0$ for all $f \in \G^1\setminus \{0\}$. 
We have written $\BB(\G^1,\G^{-1})$ for the set of bounded linear operators from 
$\G^1$ to $\G^{-1}$ and kept the notation $\langle \cdot, \cdot \rangle$
for the duality between $\G^1$ and $\G^{-1}$.
Since $S$ is positive we define the completion $\S$
of $\G^1$ with respect to the norm $\|f\|_\S := \langle f, Sf \rangle^{1/2}$.
Its adjoint space $\S^*$ can then be identified with the completion of $S\G^1$ 
with respect to the norm $\|g\|_{\S^*} := \langle g,S^{-1}g\rangle^{1/2}$. 
One observes that $S$ extends to a unitary element of $\BB(\S,\S^*)$.
$\S$ and $\S^*$ are Hilbert spaces which are generally not comparable with $\H$.
But since $\G^1 \hookrightarrow \S$ and $\S^* \hookrightarrow \G^{-1}$ it
makes sense to assume that $\{W_t\}_{t\in \R}$ restricts to a $C_0$-group in 
$\S^*$, or equivalently that it extends to a $C_0$-group in $\S$. 
Under this assumption (tacitly assumed in the sequel) we still keep the
notation $\{W_t\}_{t\in \R}$ for these $C_0$-groups.
Endowed with the graph norm, the domain of the generator of the $C_0$-group 
in $\S^*$ is denoted by $D(A,\S^*)$.
 
\begin{Definition}
For $j\in \{1,2\}$, let $\T_j$ be one of the spaces $\H, \G^s, \S$ or $\S^*$ 
introduced above.
An operator $T \in \BB(\T_1,\T_2)$ belongs to $C^1(A; \T_1, \T_2)$ if the map 
\begin{equation*}
\R \ni t \mapsto W_{-t}TW_t \in \BB(\T_1,\T_2)
\end{equation*} 
is strongly differentiable.
Its derivative at $t=0$ is denoted by $[iT,A] \in \BB(\T_1,\T_2)$.
\end{Definition}
 
Before stating the main result of this section let us recall
some known facts.
By duality and interpolation, any symmetric operator $T$ in $\H$
with $T \in \BB(\G^2, \H)$ has a unique extension 
to a symmetric element of $\BB(\G^1,\G^{-1})$, still denoted by $T$.
Then, the assumption $T \in \BB(\S,\S^*)$ has an unambiguous meaning.
It is equivalent to the requirement that $T(\G^1) \subset \S^*$
and $T: \G^1 \to \S^*$ is continuous when $\G^1$ is provided
with the topology induced by $\S$.
In that case the unique extension to a continuous mapping from $\S$ 
to $\S^*$ is still denoted by $T$.
On the other hand, if $\E$ is the Banach space 
$(D(A,\S^*), \S^*)_{1/2,1}$ defined by real interpolation 
(see for example \cite[Prop.~2.7.3]{ABG}), then one has the natural 
continuous embeddings~: 
\begin{equation*}
\BB(\G^{-1},\G^1)\subset \BB(\S^*,\S) \subset \BB(\E,\E^*)\ .
\end{equation*}

\begin{Theorem}\label{thmabstrait}
Let $H$ be a self-adjoint operator in $\H$ that belongs to $C^1(A;\G^2,\H)$
and assume that there exist two constants $c_{\hbox{\tiny \rm 1}}\geq 0$
and $c_{\hbox{\tiny \rm 2}} > 0$ such that 
\begin{equation}\label{posi}
S:=-c_{\hbox{\tiny \rm 1}}\;\!H + [iH,A]\ >\ 0
\qquad \hbox{and}\qquad [iH,A]\ \geq \ -c_{\hbox{\tiny \rm 2}}\ .
\end{equation}
Assume furthermore that $[iH,A]$ extends to an element of $C^1(A; \S,\S^*)$. 
Then, there exists $c<\infty$ such that
\begin{equation}\label{Kre}
|\langle f, (H-\lambda \mp i\mu)^{-1}f\rangle | \ \leq \ c\;\! \|f\|_{\E}^2
\end{equation} 
for all $\lambda \in \R$ with $c_{\hbox{\tiny \rm 1}}\lambda \geq 0$, 
all $\mu >0$ and all $f \in \E$.
\end{Theorem}

We observe that the condition on $\lambda$ splits into two cases.
Either $c_{\hbox{\tiny \rm 1}}=0$ and then the result holds for all 
$\lambda \in \R$, or $c_{\hbox{\tiny \rm 1}}>0$ and then $\lambda$
has to be restricted to the positive axis.
Since the case $c_{\hbox{\tiny \rm 1}}=0$ was already treated in 
\cite{BKM} we shall state two well known corollaries  
only in the case $c_{\hbox{\tiny \rm 1}}>0$.

\begin{Corollary}\label{coro}
Assume that the assumptions of Theorem \ref{thmabstrait} hold
for some $c_{\hbox{\tiny \rm 1}}>0$. Then,
\begin{itemize}
\item[{\rm (i)}] any element of $\BB\big((\E^*)^\circ, \H\big)$ is 
$H$-smooth on $[0,\infty)$, where $(\E^*)^\circ$ stands for the closure
of $\S$ in $\E^*$,
\item[{\rm (ii)}] the spectrum of $H$ on $[0,\infty)$ is absolutely
continuous.
\end{itemize}
\end{Corollary}

\section{Application to Schr\"odinger operators}

In this section, we apply the abstract result to some Schr\"odinger operators 
in the Hilbert space $\H:=L^2(\R^n)$.
Let us first recall that for $j \in \{1,\ldots, n\}$, $Q_j$ is the operator of 
multiplication by the variable $x_j$, $P_j:=-i\nabla_j$ is a component
of the momentum operator and $-\Delta \equiv P^2$ is Laplace operator
on $\R^n$. 
For each $s \in \R$, $\H^s$ denotes the usual Sobolev space of order
$s$ on $\R^n$.

Let $V$ be a real and bounded $C^\infty(\R^n)$-function. 
We shall work under this smoothness assumption that is not
essential but which simplifies our arguments.
The Schr\"odinger operator 
\begin{equation*}
H:= -\Delta + V
\end{equation*}
is self-adjoint in $\H$ with domain $\G^2 \equiv \H^2$. 
Obviously one has $\G^1 \equiv \H^1$, 
and by duality, $\G^{-2} \equiv \H^{-2}$ and $\G^{-1}\equiv \H^{-1}$.
It is well known that all these spaces are invariant under the
action of the dilation group $\{W_t\}_{t \in \R}$ whose  
generator $A$ has the form $ A:= \hbox{$\frac{1}{2}$}(P\cdot Q + Q \cdot P)$.

Let us now assume that the map $\R^n \ni x \mapsto 
\widetilde{V} (x):= \sum_{j=1}^n x_j[\partial_j V](x)\in \R$ is bounded.
It follows that $H \in C^1(A;\H^2,\H)$ and that
$[iH,A] = -2\Delta - \widetilde{V}$.
In this situation the main positivity requirement of Theorem \ref{thmabstrait}
is that there exists $c_{\hbox{\tiny \rm 1}}\geq 0 $ such that
\begin{equation*}
-(2-c_{\hbox{\tiny \rm 1}})\;\!\Delta - c_{\hbox{\tiny \rm 1}}V 
-\widetilde{V}\ >\ 0\ .
\end{equation*}
One observes that if there exists $c_{\hbox{\tiny \rm 1}}\in [0,2)$ such 
that $-c_{\hbox{\tiny \rm 1}} V - \widetilde{V}\geq 0$, then 
this inequality is obviously satisfied.

In the next proposition we use this idea and give a very simple and 
explicit application of Theorem \ref{thmabstrait}.
But let us also notice that if $n \geq 3$, some additional positivity
can be obtained from the inequality 
$-\Delta \geq (\hbox{$\frac{n-2}{2}$})^2|Q|^{-2}$.
For purposes of simplicity we do not take this improvement into account, 
and refer to \cite{BKM} for an extensive use of this inequality in the 
special case $c_{\hbox{\tiny \rm 1}}=0$.
 
\begin{Proposition}\label{application}
Let $V$ be a real and bounded $C^\infty(\R^n)$-function.
Assume furthermore that the following three conditions are satisfied 
for all $x \in \R^n$ :
{\rm (i)} $V(x) \leq 0$,
{\rm (ii)} there exists $c\in [0,2)$ such that
$|\widetilde{V}(x)| \leq -c \;\!V(x)$,
{\rm (iii)} there exists $d\geq 0$ such that~:
\begin{equation*}
\Big|\sum_{j,k=1}^n x_j \partial_j x_k \partial_k V(x)\Big| \ \leq \ 
-d\;\! V(x)\ .
\end{equation*}
Then for $c_{\hbox{\tiny \rm 1}} \in (c,2)$ fixed and
$S:= -(2-c_{\hbox{\tiny \rm 1}})\;\!\Delta - c_{\hbox{\tiny \rm 1}}\;\!V 
-\widetilde{V}$, the limiting absorption principle \eqref{Kre} is satisfied for
all $\lambda\geq 0$, all $\mu >0$ and all $f \in \E$. 
Furthermore, any element of $\BB\big((\E^*)^\circ, \H\big)$ is 
$H$-smooth on $[0,\infty)$ and the spectrum of $H$ on $[0,\infty)$ is 
absolutely continuous.
\end{Proposition}

Since the spaces $\E$, $\E^*$ and $(\E^*)^\circ$ are rather intricate, 
$H$-smooth operators are not so easily exhibited.
But under one not too restrictive extra assumption on $V$,  
a large class of $H$-smooth operators can be constructed.
For that purpose let us set 
$M(x):=\min\big\{-V(x), \hbox{$\frac{1}{|x|^2}$}\big\}$
for any $x \in \R^n$.

\begin{Corollary}\label{VVV}
Assume that $V$ satisfies the assumptions of Proposition \ref{application}
with {\rm (i)} replaced by $V(x) <0$ for all $x \in \R^n$.
If $L: \R^n \to \R$ is a Borel function that satisfies  
$|L(x)|\leq c\;\! M(x)^{\frac{1}{4}+\delta}\;\!\big(-V(x)\big)^{\frac{1}{4}-\delta}$
for some $\delta \in (0,\hbox{$\frac{1}{4}$})$, $c<\infty$ and all $x \in \R^n$, 
then the operator of multiplication by $L$ is $H$-smooth on $[0,\infty)$.
\end{Corollary}

Let us remark that if the additional assumption is replaced by the even stronger 
requirement $V(x)\leq -\e (1+x^2)^{-\mu/2}$ for some $\e>0$ and $\mu \in (0,2)$, 
then a similar result already appears in \cite[Thm.~1.9]{N} or in \cite[Cor.~3.5]{FS}.

Before starting proofs we warn the reader that the same letter 
$c$ or $d$ may denote different constants from line to line.

\begin{proof}[Proof of Proposition \ref{application}]
Let us write $\D$ for the set $C^\infty_c(\R^n)$ of smooth functions 
on $\R^n$ with compact support. 
Because of our smoothness assumption on $V$, all calculations below 
are well justified on $\D$. 

(a) One has already noticed that $H$ is a self-adjoint operator in $\H$
with domain $\H^2$, and that $H$ belongs to $C^1(A; \H^2,\H)$.
Let us now fix $c_{\hbox{\tiny \rm 1}} \in (c,2)$ and observe that~:
\begin{equation}\label{plusque}
S\ \geq \ -(2-c_{\hbox{\tiny \rm 1}})\;\!\Delta\ , \quad
S \ \geq \ - (c_{\hbox{\tiny \rm 1}}-c)\;\!V \quad \hbox{ and } \quad
S\ > \ 0\ .
\end{equation}
Furthermore the self-adjoint operator $[iH,A]=-2\Delta - \widetilde{V}$ 
is bounded from below.
Thus, both conditions in \eqref{posi} are satisfied. 

(b) By performing some easy calculations on $\D$ and by taking into account 
hypotheses (ii) and (iii) and the inequalities \eqref{plusque}
one obtains that there exists $d>0$ such that on $\D$ the 
following inequalities hold~:
\begin{eqnarray}
\label{essai1} -d\;\! S \ \leq \ [iH,A] \ \leq \ d\;\! S\ , \\
\label{essai2} -d\;\! S \ \leq \ \big[i[iH,A],A\big] \ \leq \ d\;\! S \ ,\\
\label{essai3} -d\;\! S \ \leq \ [iS,A] \ \leq \ d\;\! S\ .
\end{eqnarray}
It follows from \eqref{essai1} that  
$|\langle f,[iH,A]f\rangle |\leq d\;\!\langle f, Sf \rangle 
\equiv d\;\!\|f\|^2_\S$ for all $f \in \D$,
and then from the density of $\D$ in $\S$ that $[iH,A]$ extends
to an element of $\BB(\S,\S^*)$.
Relation \eqref{essai2} leads to the same conclusion for 
the operator $ \big[i[iH,A],A\big]$.

(c) We check now that $\{W_t\}_{t\in \R}$ extends
to a $C_0$-group in $\S$. 
This easily reduces to the proof that
$\|W_t f\|_\S \leq c(t) \|f\|_\S$ for all $f \in \D$ and $t \in \R$.
By  \eqref{essai3} one has~:
\begin{equation*}
\|W_t f\|^2_\S \ = \ \langle f, Sf\rangle + 
\int_0^t \langle W_\tau f, [iS,A]W_\tau f\rangle \;\!\de \tau
\ \leq \ \|f\|^2_\S + d \;\!\Big|\int_0^t \|W_\tau f\|_\S^2 \; \de \tau\Big|\ .
\end{equation*}
The function $(0,t) \ni \tau \mapsto \|W_\tau f\|_\S^2 \in \R$ is bounded
(since $\G^1 \hookrightarrow \S$), and hence by a simple form of the
Gronwall Lemma, we get the inequality $\|W_t f\|_\S \leq e^{\frac{d}{2}|t|}
\|f\|_\S$.
Thus $\{W_t\}_{t \in \R}$ extends to a $C_0$-group in $\S$, and by duality
$\{W_t\}_{t\in \R}$ also defines a $C_0$-group in $\S^*$. 
This finishes the proof that $[iH,A]$ extends to
an element of $C^1(A;\S,\S^*)$.
All hypotheses of Theorem \ref{thmabstrait} have been checked, 
and the statements follow from this theorem and from its corollary.
\end{proof}

\begin{proof}[Proof of Corollary \ref{VVV}]
Let $\M$ be the completion of $\D$ with respect to the norm 
$\|f\|_\M:=\|M^{-1/2}f\|$, and similarly let $\N$ be the completion 
of $\D$ with respect to the norm $\|f\|_\N:=\|(-V)^{-1/2}f\|$.
We first observe that $\M \subset \N$, 
$\M \subset D(A,\S^*)$ and $\N \subset \S^*$.
Indeed, the first continuous embedding follows directly from 
the inequality $\|f\|_\N \leq \|f\|_\M$ for all $f \in \D$.
For the second we show that
$\|f\|_{\S^*}^2 + \|Af\|_{\S^*}^2 \leq c \|f\|_\M^2$ 
for $c<\infty$ and all $f \in \D$.
From Corollary 1 of \cite{K1} and \eqref{plusque} one gets
\begin{equation}\label{bruit}
\|f\|_{\S^*}^2 \ \equiv \ \langle f, S^{-1}f\rangle \ \leq \
\hbox{$\frac{1}{c_{\hbox{\tiny \rm 1}}-c}$}\langle f, (-V)^{-1}f\rangle
\ \leq \ d\;\!\|f\|_\N^2 \ \leq \ d\;\!\|f\|_\M^2\ .
\end{equation}
Furthermore, it easily follows from \eqref{plusque} that for each 
$j \in \{1,\dots,n\}$ $P_j$ extends to an element of $\BB(\H,\S^*)$, 
and therefore~:
\begin{equation*}
\|Af\|_{\S^*} \ \leq \ c \sum_{j=1}^n \|Q_jf\| + d\;\! \|f\|_{\S^*}
\ \leq \ c'\;\!\||Q|f\| + d'\;\!\|f\|_\M \ \leq \ c'' \|f\|_\M\ .
\end{equation*}
The third embedding is also obtained from \eqref{bruit}.
One may then apply \cite[Cor.~2.6.3]{ABG}
and obtains the following relations between spaces defined by 
real interpolation~:
\begin{equation}\label{interpo}
\big(\M,\N\big)_{\theta,p} \subset \big(D(A,\S^*),\S^*\big)_{\theta,p}
\quad \forall \; \theta \in (0,1)\hbox{ and } p \in [1,\infty]\ . 
\end{equation}

In order to exhibit explicit norms on $(\M,\N)_{\theta,2}$ let us set 
$\Lambda :=\big(\hbox{$\frac{-V}{M}$}\big)^{1/2}$
and observe that $\Lambda \geq 1$. 
It is easily checked that the couple $(\M,\N)$ is quasi-linearizable in
the sense of \cite[Sec.~2.7]{ABG} (with $V_\tau := (1+\tau \Lambda)^{-1}$).
By applying then Lemma 2.7.1 of the same reference, one obtains that an 
admissible norm on $(\M,\N)_{\theta,p}$ is given by the expression 
$\big(\int_1^\infty \big\|r^{1-\theta}\hbox{$\frac{\Lambda}{r + \Lambda}$}
f\big\|^p_\N\; \hbox{$\frac{\de r}{r}$}\big)^{1/p}$.
Furthermore, by the same argument as in the proof of \cite[Prop.~2.8.1]{ABG}
one gets that in the special case $p=2$ this norm is equivalent
to the norm given by $\|\Lambda^{1-\theta}f\|_\N$.
Altogether, one has obtained that the interpolation space 
$(\M,\N)_{\theta,2}$ is equal to the completion of $\D$ with respect to
the norm $\|\Lambda^{1-\theta}(-V)^{-1/2}f\|$.

For each $\e\in (0,\hbox{$\frac{1}{2}$})$ let us set
$\theta := \hbox{$\frac{1}{2}$}-\e$ and $\F_\e:= (\M,\N)_{\theta,2}$.
One has $\F_\e\subset (\M,\N)_{1/2,1}$ 
\cite[Prop.~2.4.1]{ABG}, and it follows then by \eqref{interpo} 
that $\F_\e \subset \E$.
Thus any element of $\BB(\F_\e^*, \H)$ is $H$-smooth on $[0,\infty)$.
It is then readily checked that the operator $L$ belongs to 
$\BB(\F_\e^*, \H)$ with $\e = 2\delta$, which finishes the proof.
\end{proof}

\section{Proof of the main theorem}\label{secproof}

This section is entirely devoted to the proof of the abstract result.

\begin{proof}[Proof of Theorem \ref{thmabstrait}]
(a) For $\lambda \in \R$ with $c_{\hbox{\tiny \rm 1}} \lambda \geq 0$, 
$\mu >0$ and $\e>0$, let us consider the operators 
\begin{equation}\label{defop}
\big(H-\lambda \mp i \mu \mp i \e [iH,A]\big) \in \BB(\G^2,\H)\ .
\end{equation}
We first prove that there exists $\e_0>0$ such that for
$\e \in (0,\e_0)$ these operators are isomorphisms from $\G^2$ to $\H$.
As a consequence of the open mapping theorem, it is enough to prove
that they are bijective.
For that purpose, let us notice that for any $f \in \G^2\setminus \{0\}$ and 
$T \in \BB(\G^2,\H)$ one has $\langle f, Tf \rangle = 0$ if and only if 
$\re \langle f, Tf \rangle = 0$ and $\im \langle f, Tf \rangle = 0$.
It follows that if there exist two finite numbers $c$ and $d$ such that
$c \re \langle f, Tf \rangle + d \im \langle f, Tf \rangle \neq 0$,
then $\langle f, Tf \rangle \neq 0$.
Now, one observes that
\begin{eqnarray}\label{Bou}
\nonumber &&- c_{\hbox{\tiny \rm 1}} \re \big\langle f,
\big(H-\lambda \mp i\mu \mp i\e[iH,A]\big)f\big\rangle 
\mp \hbox{$\frac{1}{\e}$}\im \big\langle f,
\big(H-\lambda \mp i\mu \mp i\e[iH,A]\big)f\big\rangle \\
&=&  \big\langle f, \big(-c_{\hbox{\tiny \rm 1}}  
H + [iH,A] \big)f\big\rangle + 
( c_{\hbox{\tiny \rm 1}} \lambda + \hbox{$\frac{\mu}{\e}$})
\|f\|^2 \ \geq \ \langle f,S f\rangle \ >\ 0\ ,
\end{eqnarray}
which implies that $\langle f,(H-\lambda \mp i \mu \mp i \e [iH,A]) f\rangle
\neq 0$ for all $f \in \G^2 \setminus \{0\}$.
Thus both operators in \eqref{defop} are injective.
Furthermore, for $\e \in (0,\e_0)$ with $\e_0 := (\|[iH,A]\|_{\G^2 \to \H})^{-1}$
one easily deduces that they are closed operators in $\H$
and adjoint to each other.
This immediately leads to their surjectivity \cite[Sec.~V.3.1]{K2}
and thus to their bijectivity.

(b) For $\e \in (0,\e_0)$ let us set 
$G_\e^\pm := (H-\lambda \mp i\mu \mp i \e [iH,A])^{-1}$.
These operators belong to $\BB(\H,\G^2)$, and by duality and interpolation 
to $\BB(\G^{-1},\G^1)\subset \BB(\S^*,\S)$.
It is then easily shown that for all $f,g \in \G^{-1}$ :
$\langle f,G_\e^\pm g\rangle \ = \ \langle G_\e^\mp f,g\rangle$.
By taking into account these equalities and the continuous extensions of the 
inequalities \eqref{Bou} valid for all $f \in \G^1 \setminus \{0\}$, one 
observes that there exists $c<\infty$ such that for all $f \in \G^{-1}$~:
\begin{eqnarray*}
&& \|G_\e^\pm f\|^2_\S \  = \ \langle G_\e^\pm f, S G_\e^\pm f \rangle \\
& \leq & c_{\hbox{\tiny \rm 1}} \big| \re \big \langle G_\e^\pm f, 
\big(H-\lambda \pm i\mu \pm i\e[iH,A]\big) G_\e^\pm f \big\rangle \big| \\
&& + \; \hbox{$\frac{1}{\e}$} 
\big| \im \big \langle G_\e^\pm f, 
\big(H-\lambda \pm i\mu \pm i\e[iH,A]\big) G_\e^\pm f \big\rangle \big| \\
& \leq & c_{\hbox{\tiny \rm 1}} |\langle f,G_\e^\pm f\rangle| + 
\hbox{$\frac{1}{\e}$} |\langle f,G_\e^\pm f\rangle|
\ \leq \ \hbox{$\frac{c}{\e}$}|\langle f,G_\e^\pm f\rangle|\ .
\end{eqnarray*}
Thus, one has obtained that for all $f \in \G^{-1}$~:
\begin{equation}\label{f1}
\|G_\e^\pm f \|_\S \ \leq \ \hbox{$\sqrt{\frac{c}{\e}}$}\ |\langle f, 
G_\e^\pm f \rangle |^{1/2} \ ,
\end{equation}
and by using the inequality 
$|\langle f, g\rangle| \leq \|f\|_{\S^*} \|g\|_\S$ valid for
all $f \in \S^*$ and $g \in \S$, it follows that~:
\begin{equation}\label{f2}
\|G_\e^\pm\|_{\S^* \to \S} \ \leq \ \hbox{$\frac{c}{\e}$} \ .
\end{equation}

(c) This part of the proof is similar to parts (ii) to (iv) of the proof of 
\cite[Thm.~2.1]{BM}.
For $\e >0$ and $f\in \E \equiv (D(A,\S^*),\S^*)_{1/2,1}$, 
let us set $f_\e:= \hbox{$\frac{1}{\e}$}\int_0^\e (W_t f) \;\!\de t$.
Then, $f_\e \in D(A,\S^*)$, $\e \mapsto f_\e \in \S^*$ 
is $C^1$ in norm, $f_\e \to f$ in $\S^*$ as $\e \to 0$ and 
\begin{equation}\label{int}
\int_0^1 \big(\|f'_\e\|_{\S^*} + \|Af_\e\|_{\S^*}\big)\;\e^{-1/2}\;\de \e 
\ \leq \ c \|f\|^2_\E\ ,
\end{equation} 
where $f'_\e$ denotes the derivative of the map $\e \mapsto f_\e$
(see \cite[Thm.~2.1]{BM}).

Now, for $\e \in (0,\e_1)$, with $\e_1 := \min\{\e_0,1\}$, let us set
$F_\e^\pm := \langle f_\e, G_\e^\pm f_\e \rangle$. 
A formal calculation leads to~:
\begin{equation*}
\hbox{$\frac{\de}{\de \e}$}F_\e^\pm \ \equiv \ (F_\e^\pm)' \ = \ 
\langle f'_\e \mp Af_\e, G_\e^\pm f_\e\rangle + 
\langle G_\e^\mp f_\e, f'_\e \pm Af_\e\rangle
-\e \langle G_\e^\mp f_\e, \big[i[iH,A],A\big] G_\e^\pm f_\e\rangle\ .
\end{equation*}
A rigorous proof of these equalities can be derived similarly as in the 
usual Mourre theory, see for example \cite[Sec.~7.3]{ABG} or 
\cite[Lem.~3.4]{BGM}. By taking \eqref{f1} into account we obtain the 
fundamental differential inequalities~:
\begin{equation}\label{eqineq}
\hbox{$\frac{1}{c}$}|(F_\e^\pm)'| \ \leq \ \hbox{$\frac{1}{\sqrt{\e}}$}
\big(\|f'_\e\|_{\S^*} + \|Af_\e\|_{\S^*}\big)|F_\e^\pm|^{1/2} +
\big\| \big[i[iH,A],A\big] \big\|_{\S\to\S^*}|F_\e^\pm|\ .
\end{equation}

Then, by an application of the Gronwall lemma \cite[Lem.~7.A.1]{ABG} 
together with the use of the inequality \eqref{int}
one concludes that $F_0^\pm = \lim_{\e \to 0}F_\e^\pm$ exist and satisfy 
$|F_0^\pm| \leq c (|F_{\e_1}^\pm| + \|f\|^2_\E)$.
Furthermore, one has by \eqref{f2} that
\begin{equation*}
|F_{\e_1}^\pm| \ \leq \ \|G_{\e_1}^\pm\|_{\S^* \to \S}\;\!\|f_{\e_1}\|^2_{\S^*} 
\ \leq \ c \Big[\int_0^{\e_1} \|W_t f\|_{\S^*}\;\!\de t\Big]^2
\ \leq \ d \|f\|^2_{\S^*} \ \leq \ d \|f\|^2_\E\ ,
\end{equation*} 
which leads to the expected inequalities :
$|F_0^\pm| \ \leq \ c \|f\|^2_\E$.

(d) It only remains to show that
$F_0^\pm = \langle f, (H-\lambda \mp i\mu)^{-1} f\rangle$, 
{\it i.e.}~that the right objects have been obtained.
For that purpose, let us set $G_0^\pm:=(H-\lambda \mp i\mu)^{-1}$ and 
observe that
\begin{eqnarray*}
|\langle f_\e, G_\e^\pm f_\e\rangle - \langle f, G_0^\pm f\rangle|
& \leq & \|f_\e - f\|_{\S^*} \;\|G_\e^\pm \|_{\S^* \to \S} \;\|f_\e\|_{\S^*} \\
&& + \;\|f\|_{\S^*}\; \|G_\e^\pm -G_0^\pm \|_{\S^* \to \S}\;\|f_\e\|_{\S^*} \\
&&+\; \|f\|_{\S^*}\;\|G_0^\pm \|_{\S^*\to\S}\;\|f_\e - f\|_{\S^*}\ .
\end{eqnarray*}
Since $\|f_\e-f\|_{\S^*}\to 0$ as $\e \to 0$ and $\|T\|_{\S^*\to\S} \leq
d\;\!\|T\|_{\G^{-1}\to\G^1}$ for all $T \in \BB(\G^{-1},\G^1)$, 
it is enough to prove that
for $\lambda$ and $\mu$ fixed there exist $\e_2>0$ and $c<\infty$ 
such that $\|G_\e^\pm \|_{\G^{-1} \to \G^1} \leq c$ for all $\e \in [0,\e_2]$. 
Indeed, by using the second identity of the resolvent 
one then gets the inequalities~:
\begin{equation*}
\|G_\e^\pm -G_0^\pm \|_{\G^{-1}\to\G^1} \ \leq \  
\|G_\e^\pm \|_{\G^{-1}\to \G^1} \;\!\|\e [iH,A]\|_{\G^1 \to \G^{-1}} \;\!
\|G_0^\pm \|_{\G^{-1}\to\G^1} \ \leq \ \e\;\! c^2 \;\!
\|[iH,A]\|_{\G^1 \to \G^{-1}}\ .
\end{equation*}

So, let us set 
$\e_2:= \min\{\e_1,\hbox{$\frac{\mu}{2c_{\hbox{\tiny \rm 2}}}$}\}$
and observe that for $\e \in [0,\e_2]$ and all $f \in \G^2$, the inequality
$\hbox{$\frac{\mu}{2}$} \|f\|^2 + \langle f, \e[iH,A]f\rangle \geq 0$ holds.
It easily follows that for all $f \in \G^2$~:
\begin{equation*}
\big\|\big(H-\lambda \mp i\hbox{$\frac{\mu}{2}$} 
\mp i(\hbox{$\frac{\mu}{2}$}+\e[iH,A])\big)f\big\| \ \geq \ 
\big\|\big(H-\lambda \mp i\hbox{$\frac{\mu}{2}$}\big)f\big\|\ ,
\end{equation*}
and then that for all $g \in \H$~:
\begin{equation*}
\big\|(H-\lambda \mp i\hbox{$\frac{\mu}{2}$})G_\e^\pm g\big\|\ \equiv \ 
\big\|(H-\lambda \mp i\hbox{$\frac{\mu}{2}$})
\big(H-\lambda \mp i\hbox{$\frac{\mu}{2}$} 
\mp i(\hbox{$\frac{\mu}{2}$}+\e[iH,A])\big)^{-1}g\big\| \leq \|g\|\ .
\end{equation*}
Therefore, one obtains that for all $g \in \H$~:
\begin{eqnarray*}
&&\big\|G_\e^\pm g\big\|_{\G^2}
\ = \  \big\|(H+i)G_\e^\pm g\big\| \ = \ 
\big\|(H+i)(H-\lambda \mp i\hbox{$\frac{\mu}{2}$})^{-1} 
(H-\lambda \mp i\hbox{$\frac{\mu}{2}$})G_\e^\pm g \big\| \\
&\leq& \big\|(H+i)(H-\lambda \mp i\hbox{$\frac{\mu}{2}$})^{-1}\big\|\;\!\|g\|
\ \leq\  \big(1+\big|\lambda \pm i \hbox{$\frac{\mu}{2}$} +i\big|
\hbox{$\frac{2}{\mu}$}\big)\;\!\|g\|\ ,
\end{eqnarray*}
or equivalently that $\|G_\e^\pm\|_{\H \to \G^2}
\leq c$ with $c$ independent of $\e \in [0,\e_2]$.
By duality and interpolation, one concludes that
$\|G_\e^\pm\|_{\G^{-1}\to \G^1} \leq c$ for all $\e \in [0,\e_2]$.
\end{proof}

\section*{Acknowledgements}
The author thanks the Swiss National Science Foundation for
its financial support.

\end{document}